\documentclass[12pt]{article}
\usepackage{latexsym}
\newtheorem{Lemma}{Lemma}
\newtheorem{Theorem}{Theorem}
\newtheorem{Prop}{Proposition}

\newtheorem{Cor}{Corollary}

\newcommand{\pf}{\medskip\noindent{\sc Proof: }}
\newcommand{\qed}{$\Box$}
\newcommand{\rTo}{\,\longrightarrow\,}

\newcommand{\DS}{\displaystyle}
\begin{document}
\title{Factorials and powers, a minimality result}
\author{David E. Radford \\
Department of Mathematics, Statistics \\
and Computer Science (m/c 249)    \\
851 South Morgan Street   \\
University of Illinois at Chicago\\
Chicago, Illinois 60607-7045
 }
\maketitle
\date{}

\begin{abstract}
{\small  \rm Let $a > 1$. Then $a^n < n!$ for some positive integer $n$. We show that the smallest such $n$ is one of a pair of possibilities, or is one possibility, which we show how to calculate. There are three interesting numerical sequences which play a central role in our arguments. This paper is based on the improvement on Sterling's approximation of factorials due to Robbins \cite{Robbins} and results of \cite{Radford}.}
\end{abstract}
\let\thefootnote\relax\nonumber
%SECTION
%
\section{Introduction}\label{SecIntroduction}
Let $a > 1$. Then $a^n \leq n!$ for some $n \geq 1$. Let $n_a$ be the smallest such positive integer. Then $n_a$ is the unique solution $n$ to $n \geq 2$ and $a$ belongs to the half open interval $(\sqrt[n-1]{(n-1)!}\,, \; \sqrt[n]{n!}]$. Now $a$ also belongs to a half open interval $(n/e, \; (n+1)/e]$ for a unique $n \geq 2$, an interval whose length is $1/e$. We describe $n_a$ in terms of $n$. Let $\sigma_n$ be the smallest integer such that $e\sqrt[n]{n!} < n + \sigma_n$. Then $\sigma_n \geq 2$ for all $n \geq 1$ and $\sigma_n < n$ when $n \geq 3$. By Theorem \ref{TheoremMain}, the main result of this paper, $n_a = n-m+1, n-m+2$, or $n-m+3$, where $n \geq 3$ and $m = \sigma_n$. More precisely, the theorem shows how the segment $\sigma_{n-m}, \ldots, \sigma_n$ determines pairs of possibilities for $n_a$ and in one case $n_a$ precisely.

The sequence $\sigma_1, \sigma_2, \sigma_3, \ldots$ is increasing and each term is equal to its successor or is one less. The sequence appears to grow very, very slowly. For example $\sigma_n = 15$ when $n = 10^{12}$. A key ingredient of the proof of Theorem \ref{TheoremMain} is Corollary \ref{CorOneOrTwoValuesB} which states that when $n \geq 3$ that the sequence $\sigma_{n-m}, \ldots, \sigma_n$ has one or two values, where $m = \sigma_n$. There are a number of interesting sequences associated with analysis of $n_a$ which we examine in this paper.

We explore in detail functions related to an approximation of $e\sqrt[n]{n!}$ which is equivalent to an improvement of Sterling's approximation of factorials due to Robbins \cite{Robbins}. Consider the four functions
$$
L(x) = \DS{\pi^\frac{1}{2x}}e^{\frac{1}{(12x+1)x}}, \;\; R(x) = \pi^\frac{1}{2x}e^\frac{1}{12x^2}, \;\; P(x) = x^{\frac{1}{x}}, \;\; \mbox{and} \;\; T(x) = xP(x)
$$
which are defined and differentiable on $(0, \infty)$. Our way of writing Robbins's improvement is $(1/2)L(n)T(2n) < e\sqrt[n]{n!} < (1/2)R(n)T(2n)$ for all $n \geq 1$.

The study of the four functions mentioned above is based on elementary estimates of $a^{1/c}$, where $a > 1$ and $c > \ln a$. These functions, as well as the derivatives $T'(x)$ and $T''(x)$, and estimates are discussed in Section \ref{SecERootNFactorial}. The function $T'(x)$ is of particular importance to us.

Let $F(x)$ be any one of the four functions $L(x)$, $R(x)$, $P(x)$, and $T'(x)$. Since $\lim_{x \rTo \infty}F(x) = 1$, we find it useful to write $\DS{F(x) = 1 + \frac{a_F(x)}{x}}$, where $a_F(x) = (F(x)-1)x$, and develop properties of $a_F(x)$. This is done in Section \ref{SecBasicRevisited} with particular emphasis on $L(x)$ and $R(x)$.

In the very short Section \ref{SecTDifffences} we consider our first sequence associated with the analysis of $n_a$ which is $s_1, s_2, s_3, \ldots \,$, where
$$
s_n = (n+1)\sqrt[n+1]{n+1} - n\sqrt[n]{n} = T(n+1) - T(n) = T'(x)
$$
for some $n < x < n+1$ by the Mean Value Theorem. Set $\textsl{a}_s(x) = (T'(x)-1)x$ for $x > 0$. Theorem \ref{TheoremNRootN}, the proof of which is based on our analysis of $T'(x)$, lists properties of the sequence $s_1, s_2, s_3, \ldots \,$, in particular $2 > s_1 > s_2 > s_3 > \cdots > 1$,
$$
1 + \frac{\textsl{a}_s(n+1)}{n+1} < s_n < 1 + \frac{\textsl{a}_s(n)}{n}
$$
for all $n \geq 1$, and $\lim_{n \longrightarrow \infty}\textsl{a}_s(n) = 1$.

Let $x > 0$ and let $G(x) = \Gamma(x+1)^{1/x}$, where $\Gamma(x)$ is the gamma function. In Section \ref{SecTauDifffences}  we analyze the second sequence in connection with $n_a$ which is $S_1, S_2, S_3, \ldots \,$, where
$$
S_n = e\sqrt[n+1]{(n+1)!} - e\sqrt[n]{n!} = eG(n+1) - eG(n) = eG'(x)
$$
for some $n < x < n+1$. Set $\textsl{a}_S(x) = (eG'(x) - 1)x$. This short section is about \cite[Theorem 3]{Radford} which states in particular $1.15 > S_1 > S_2 > S_3 > \cdots > 1$,
$$
1 + \frac{\textsl{a}_S(n+1)}{n+1} < S_n < 1 + \frac{\textsl{a}_S(n)}{n}
$$
for all $n \geq 18$, and $\DS{\lim_{n \longrightarrow \infty}\textsl{a}_S(n) = \frac{1}{2}}$. We use the sequence of Section \ref{SecTDifffences} to establish an important detail for $S_1, S_2, S_3, \ldots \;$. Note that \cite[Theorem 3]{Radford} has implications for the half open intervals mentioned in the fourth line of the introduction, namely that their lengths exceed $1/e$ which they approach as $n$ becomes large.

The sequences of Sections \ref{SecTDifffences} and \ref{SecTauDifffences} are of a more general type which we study in Section \ref{SecsigmaSequence}. Let $T_1, T_2, T_3, \ldots \;$ be any sequence of real numbers which satisfy $T_1 \geq 1$, $1 \leq s_n < 2$ for all $n \geq 1$, where $s_n = T_{n+1} - T_n$, and there is an $a > 0$ such that $s_n < 1 + a/n$ for all $n \geq 1$. The sequence of Section \ref{SecTDifffences} fits into this framework, where $T_n =  n\sqrt[n]{n}$, and the sequence of Section \ref{SecTauDifffences} does as well, where $T_n = e\sqrt[n]{n!} \;$.

It is easy to see that $n \leq T_n < T_1 + 2(n-1)$ for all $n \geq 1$. Let $\sigma_n$  be the largest integer $\ell$ such that $n + (\ell -1) \leq T_n$, or equivalently the smallest integer $\ell$ such that $T_n < n + \ell$. Then $1 \leq \sigma_1 \leq \sigma_2 \leq \sigma_3 \leq \cdots \;$ and $\sigma_n < \sigma_{n+1}$ implies $\sigma_{n+1} = \sigma_n + 1$. The proof of the main result of this paper, again Theorem \ref{TheoremMain}, is based on properties of the sequence $\sigma_1, \sigma_2, \sigma_3, \ldots \;$, our third and last sequence associated with the analysis of $n_a$.

Now let $T_1, T_2, T_3, \ldots \;$ be as defined for the sequence of Section \ref{SecTDifffences} or of Section \ref{SecTauDifffences}. Suppose $n \geq 3$. Then $2 \leq m = \sigma_n < n$. If the segment $\sigma_{n-m}, \ldots, \sigma_n$ has one value then
$$
n-1 \leq T_{n-m} < n \leq T_{n-m+1} < n+1 \leq T_{n-m+2} < n+2
$$
by Lemma \ref{LemmaSigmasSame} which relates the two types of half open intervals described in the first paragraph of the introduction. The segment has one or two values by Corollaries \ref{CorOneOrTwoValuesA} and \ref{CorOneOrTwoValuesB}. When the segment has two values, Proposition \ref{PropSigmasTwoValues} lists three possible outcomes similar to that of Lemma \ref{LemmaSigmasSame} described above. Lemma \ref{LemmaSigmasSame} and Proposition \ref{PropSigmasTwoValues} provide the basis for the proof of Theorem \ref{TheoremMain}. In Section \ref{SecFactorialsandPowers} we make some general observation on factorials and powers and prove Theorem \ref{TheoremMain}.

This paper was inspired by a question about the gamma function \cite{Artin} posed by a friend \cite{Wells} resulting from an investigation by his young son into the connection between factorials and powers. Numerical computations in this paper were made using the software MATLAB \cite{MATLAB}.
%
%SECTION
%
\section{Basic functions related to $e\sqrt[n]{n!}$ estimates}\label{SecERootNFactorial}
We begin with the improvement on Sterling's approximation of factorials due to Robbins \cite{Robbins} which is
\begin{equation}\label{EqRobbins}
\sqrt[]{2\pi} n^{n+ \frac{1}{2}}e^{-n}e^{\frac{1}{12n+1}} < n! < \sqrt[]{2\pi} n^{n+ \frac{1}{2}}e^{-n}e^{\frac{1}{12n}}
\end{equation}
for $n \geq 1$. Let $x > 0$. We define four functions
$$
L(x) = \DS{\pi^\frac{1}{2x}}e^{\frac{1}{(12x+1)x}}, \;\; R(x) = \pi^\frac{1}{2x}e^\frac{1}{12x^2}, \;\; P(x) = x^{\frac{1}{x}}, \;\; \mbox{and} \;\; T(x) = xP(x).
$$
An equivalent formulation of (\ref{EqRobbins}) is
\begin{equation}\label{EqCapSN}
(1/2)L(n)T(2n) < e\sqrt[n]{n!} < (1/2)R(n)T(2n)
\end{equation}
for all $n \geq 1$.

We examine these four functions in some detail. First of all observe that $\lim_{x \longrightarrow \infty} L(x) = 1 =  \lim_{x \longrightarrow \infty} R(x)$. It is easy to see that $1 < L(x) < R(x)$.
\begin{equation}\label{EqLnRn}
R(y) < L(x) \;\; \mbox{for all $y \geq x + 1/24$.}
\end{equation}
To verify this inequality, we need only show $\pi^{\frac{1}{2y} - \frac{1}{2x}} < e^{\frac{1}{12x^2+x} - \frac{1}{12y^2}}$ which has the form $\pi^a < e^b$, where $a < 0 < b$, under our conditions. Therefore $R(y) < L(x)$ as $\pi > 1$ and thus $a(\ln \pi) < 0 < b$.

Set
$$
\ell(x) = \DS{\frac{\ln \pi}{2x^2} + \frac{24x+1}{(12x + 1)^2x^2}} \;\; \mbox{and} \;\; r(x) = \DS{\frac{\ln \pi}{2x^2} + \frac{1}{6x^3}}.
$$
Writing the functions $L(x)$ and $R(x)$ in the form $e^{c(x)}$ one can see that
\begin{equation}\label{EqLDerivative}
L'(x) = -\ell(x)L(x) \;\; \mbox{and} \;\; R'(x) = -r(x)R(x).
\end{equation}
At this point it is easy to see that $L(x)$, $\ell(x)$, $R(x)$, and $r(x)$ are positive valued strictly decreasing functions on $(0, \infty)$, $\ell(x) < r(x)$, and
\begin{equation}\label{EqSmallLnSmallRn}
r(y) < \ell(x) \;\; \mbox{for all $y > x$.}
\end{equation}
Since $P(x) = x^{\frac{1}{x}} = e^{\frac{\ln x}{x}}$ it is easy to see that $P(x) > 0$ and when $x > 1$ that $P(x) > 1$. Now
%\begin{equation}\label{EqLHospitalLNXA}
$\DS{\lim_{x \longrightarrow \infty} \frac{\ln x}{x^a} = 0}$
%\end{equation}
for $a > 0$ by L'H\^{o}pital's rule for the indeterminant form $\DS{\frac{\infty}{\infty}}$. As a result $\lim_{x \longrightarrow \infty} P(x) = 1$. Next observe that $P'(x) = p(x)P(x)$, where $p(x) =  \DS{\frac{1 - \ln x}{x^2}}$. Consequently $P(x)$ is strictly increasing on $(0, e]$ and is strictly decreasing on $[e, \infty)$.

Note $T(x) > 0$ and when $x > 1$ that $T(x) > x$. It is easy to see that
\begin{equation}\label{EqTPrime}
T'(x) = \DS{P(x)\left(\frac{x+1 - \ln x}{x}\right) = T(x)\left(\frac{x+1 - \ln x}{x^2}\right)}.
\end{equation}
Now $x + 1 - \ln x > 1$ as $\ln x < x$. When $x > 1$ we note that $x + 1 - \ln x > 2$ since $\ln x < x - 1$ in this case. In particular $T(x)$ is a strictly increasing function on $(0, \infty)$. Since
\begin{equation}\label{EqTPrimePrime}
T''(x) = \DS{T(x)\left(\frac{(\ln x - 1)^2 - x}{x^4} \right)}
\end{equation}
it follows that $T''(x) < 0$ for all $x > 1$ as $(\ln x - 1)^2 < x$ in this case. The inequality is easily seen to hold when $1 < x \leq e$. Suppose $x > e$. Then $\ln x - 1 > 0$ and thus the inequality we wish to establish will follow from $\ln x < 1 + a\sqrt{x}$ for $x > 0$, where $a = 0.447$. Let $a > 0$. By computing the derivative of the difference $d(x) = 1 + a\sqrt{x} - \ln x$, where $x > 0$, we see that $d(x)$ has an absolute minimum at $x = 4/a^2$ and $d(4/a^2) > 0$ if and only if $\ln(4/a^2) < 3$, or equivalently $a > 2/e^\frac{3}{2}$, in which case $d(x) > 0$ for all $x > 0$. As $0.447 > 2/e^\frac{3}{2}$, we have shown $\ln x < 1 + 0.447\sqrt{x}$ for $x > 0$ and indeed $(\ln x - 1)^2 < 0.2x$ for $x \geq e$.

We have shown that $T''(x) < 0$ for all $ x > 1$. Therefore $T'(x)$ is strictly decreasing on $[1, \infty)$. Since $\lim_{x \longrightarrow \infty}T'(x) = 1$, we conclude that $T'(x) > 1$ for all $x \geq 1$.
\begin{Lemma}\label{LemmaRNOverLN}
Let $a > 1$ and $c > \ln a$. Then $\DS{\frac{c + \ln a}{c} < a^{1/c} < \frac{c}{c - \ln a}}$. Hence $\DS{1 + \frac{\ln a}{c} < a^{1/c} < 1 + \frac{\ln a}{c - \ln a}}$ and $\DS{1 - \frac{\ln a}{c} < a^{-1/c} < 1 - \frac{\ln a}{c + \ln a}}$.
\end{Lemma}

\pf We need only establish the first assertion.  Consider the special case $a = e$. Here $c > \ln e = 1$ and the inequalities are $\DS{1 + \frac{1}{c} < e^{1/c} < \frac{c}{c-1}}$. The first follows from the series expansion of $e^x$ and the equivalent of the second, $\DS{\frac{1}{c} < \ln (1 + \frac{1}{c-1})}$, is $\DS{\frac{x}{1+x} < \ln (1 + x)}$, where $x > 0$, with $\DS{x = \frac{1}{c-1}}$.

We have established the inequalities in the special case $a = e$. The general case follows from the special case as $a^{1/c} = e^{(1/c)\ln a} = e^{1/c'}$, where $c' = \DS{\frac{c}{\ln a} > 1}$.
\qed
\medskip

Lemma \ref{LemmaRNOverLN} implies $\DS{1 + \frac{\ln x}{x} < x^{1/x} < \frac{x}{x - \ln x}}$ for $x > 1$ which provides estimates of $P(x) = x^{1/x}$ for such $x$. There are improvements when $x$ is large enough. More precisely:
\begin{Prop}\label{PropXToOneOverX}
$\DS{1 + \frac{\ln x}{x-1} < x^{1/x}}$ when $x \geq 8.0845$ and $\DS{x^{1/x} < \frac{x+1}{x+1 - \ln x}}$ when $x \geq 6.7537$.
\end{Prop}

\pf
Let $x > e$. Let $f(x)$ be a continuous function on $(e, \infty)$. The function $H(x) = (\ln x)^2 - f(x)(\ln x) + 1$ will play a central role in our proof in two specific cases.

Now let $y = \ln x$. Then $y > 1$ and $x = e^y$. Set $\DS{K(y) = e^y\left(1 - \frac{y}{(y-1)^2}\right)}$. Then $K(y)$ is a strictly increasing function on its domain $(1, \infty)$, since its derivative is always positive, and its range is the set of real numbers.

We first establish the second inequality
\begin{equation}\label{EqNOneOverNEst}
x^{1/x} < \frac{x+1}{x+1-\ln x}
\end{equation}
for $x \geq 6.7537$. Now (\ref{EqNOneOverNEst}) holds if and only if $e^{(\ln x)/x} < e^{\ln ((x+1)/(x+1 - \ln x))}$ which holds if and only if
$$
G(x) = \DS{\frac{\ln x}{x} - \ln \left(\frac{x+1}{x+1-\ln x}\right) < 0}.
$$
Observe that
$$
G'(x) = \DS{\frac{H(x)}{x^2(x + 1 - \ln x)}},
$$
where $f(x) = 3 - 1/(x+1)$, and
$$
H(x) = H(e^y) = \DS{\frac{e^y(y^2-3y+1) + (y-1)^2}{e^y + 1} = \frac{(K(y)+1)(y-1)^2}{e^y+1}}.
$$
Therefore $G'(x) = 0$ (respectively ($G'(x) < 0$, $G'(x) > 0$) if and only if $K(y)= -1$ (respectively $K(y) < -1$, $K(y) > -1$). As a result there is a unique $x_0$ which satisfies $x_0 > e$ and $G'(x_0) = 0$. Furthermore $G'(x) > 0$ for $x > x_0$ and $G'(x) < 0$ for $e < x < x_0$. Therefore $G(x)$ is strictly increasing on $[x_0, \infty)$ and is strictly decreasing on $(e, x_0]$. Since $\lim_{x\longrightarrow \infty} G(x) = 0$ it follows that $G(x) < 0$ for all $x \geq x_0$. Since $G(6.7537) < 0$, it must be the case that $G(x) < 0$ for all $x \geq 6.7537$. We have established (\ref{EqNOneOverNEst}) for $x \geq 6.7537$.

Note that the inequality fails for $x = 6.7536$ as $G(6.7536) > 0$. It might be of interest to know that $12.5690 < x_0 < 12.5691$ which follows from the calculation $K(\ln 12.5690) < -1 < K(\ln 12.5691)$.

We now establish the first inequality
\begin{equation}\label{EqNOneOverNEstLeft}
1 + \frac{\ln x}{x - 1} < x^{1/x}
\end{equation}
for $x \geq 8.0845$. Our argument is modelled on the proof for the second inequality. We provide a rough sketch of the details.

First note that (\ref{EqNOneOverNEstLeft}) holds if and only
$$
G(x) = \DS{\frac{\ln x}{x} - \ln \left(1 + \frac{\ln x}{x - 1}\right) > 0}.
$$
Now
$$
G'(x) = \DS{-\frac{H(x)}{x^2(x - 1 + \ln x)}},
$$
where $f(x) = 3 + 1/(x-1)$.

Observe that $G'(x) = 0$ (respectively $G'(x) < 0$, $G'(x) > 0$) if and only if $K(y) = 1$ (respectively $K(y) > 1$, $K(y) < 1$). Therefore there is a unique $x_0$ which satisfies $e < x_0$ and $G'(x_0) = 0$. Also $G'(x)$ is strictly decreasing on $[x_0, \infty)$ and is strictly increasing on $(e, x_0]$. As $\lim_{x\longrightarrow \infty}, G(x) = 0$ necessarily $G(x) > 0$ for $x \geq x_0$. Since $G(8.0845) > 0$ we have established (\ref{EqNOneOverNEstLeft}) for $x \geq 8.0845$.

As $G(8.0844) < 0$ this inequality fails for $x = 8.0844$. The reader might find it interesting to know that $14.9063 < x_0 < 14.9064$ which follows from the calculation $K(\ln 14.9063 ) < 1 < K(\ln 14.9064)$.
\qed
%
%SECTION
%
\section{Basic functions revisited}\label{SecBasicRevisited}
Let $a > 0$ and $F(x)$ be a function defined on $(a, \infty)$. Let $x > a$ and set $a_F(x) = (F(x) - 1)x$. Then $\DS{F(x) = 1 + \frac{a_F(x)}{x}}$. If $F(x)$ is differentiable on $(a, \infty)$ then so in $a_F(x)$.

Now suppose $f(x)$ is a positive valued differentiable function on $(a, \infty)$ and set $F(x) = e^{f(x)}$. Then $F(x) > 1$ and $F'(x) = f'(x)F(x)$. Observe that
\begin{equation}\label{EqAFFPrime}
a_F'(x) = (xf'(x)+1)F(x) - 1.
\end{equation}

The four functions $L(x)$, $R(x)$, $P(x)$, and $T'(x)$, defined in Section \ref{SecERootNFactorial} are examples of such an $F(x)$ described in the preceding paragraph. In this section we study $a_F(x)$ and its derivative for the first two functions of these functions in detail.

We consider $L(x)$ and $R(x)$ together. We will show for them that
\begin{equation}\label{EqXFPPOne}
xf'(x) + 1 > 0
\end{equation}
and $f''(x)$ exists for all $x > a$, where $a$ is to be determined. Set
$$
G(x) = \ln ((xf'(x) + 1)F(x)) = \ln (xf'(x) + 1) + f(x).
$$
Then $G(x) = 0$ (respectively $G(x) > 0$, $G(x) < 0$) if and only if $a'_F(x) = 0$ (respectively $a'_F(x) > 0$, $a'_F(x) < 0$). Observe that
\begin{equation}\label{EqGPrimeGeneral}
G'(x) = \DS{\frac{xf''(x) + f'(x)(xf'(x) + 2)}{xf'(x) + 1}}
\end{equation}
for all $x > a$.

We note that $L(x)$ and $R(x)$ are products whose factors have the form $F(x)$, where $\DS{f(x) = \frac{A}{Bx^2 + Cx}}$ for $x > a$, and $A$, $B$, and $C$ are non-negative constants such that $A$ and one or both of $B$ and $C$ are positive. Let $f(x)$ be such a function, where $x > 0$. Then $f(x) > 0$ and $\lim_{x \longrightarrow \infty} f(x) = 0$. Thus $F(x) = e^{f(x)} > 1$ for all $x > 1$ and $\lim_{x \longrightarrow \infty} F(x) = 1$.

Let $g(x) = \DS{- \left(\frac{2Bx + C}{Bx^2 + Cx}\right)}$. Then $f'(x) = g(x)f(x)$. It is easy to see that $\lim_{x \longrightarrow \infty}xf'(x) = 0$ also and (\ref{EqXFPPOne}) holds if and only if
\begin{equation}\label{EqXFPPOneHolds}
A\left(\frac{2Bx^2 + Cx}{(Bx^2 + Cx)^2}\right) < 1.
\end{equation}
The coefficient of $A$ in (\ref{EqXFPPOneHolds}) is a strictly decreasing function which approaches $0$ as $x$ becomes large. Choose $a > 0$ such that (\ref{EqXFPPOneHolds}) holds for $x > a$. Note that
$$
xf''(x) + 2f'(x) = (x(g'(x) + g(x)^2) + 2g(x))f(x) = \DS{\left(\frac{2B^2x^3}{(Bx^2+Cx)^2}\right)f(x)}
$$
for $x > 0$. As $f(x) > 0 > f'(x)$ for $x > a$, it follows $G'(x) > 0$ for $x > a$ which means that $G(x)$ is strictly increasing on $(a, \infty)$. Since $\lim_{x \longrightarrow \infty} G(x) = 0$ it follows that $G(x) < 0$ for all $x > a$. We have shown that $a_{F}(x)$ is a strictly decreasing function on $(a, \infty)$.

Let $c(x) = (Bx^2 + Cx)/A$. Then $F(x) = e^{1/c(x)}$; for all $x > a$ it follows that $c(x) > 1$ by (\ref{EqXFPPOneHolds}). Therefore
\begin{equation}\label{EqAFEstimates}
\DS{\frac{A}{Bx+C} < a_F(x) < \frac{A}{Bx+C-A/x}}
\end{equation}
for all $x > a$ by Lemma \ref{LemmaRNOverLN}. As a result
\begin{equation}\label{EqLimitAF}
\lim_{x \longrightarrow \infty}a_F(x) = \left\{\begin{array}{lll} 0 & : & B \neq 0 \\ A/C  & : & B = 0 \end{array} \right. \,.
\end{equation}
We consider three special cases. Assertions are easy to verify.

$F(x) = \pi^\frac{1}{2x}$. Here $A = \ln \pi$, $B = 0$, and $C = 2$. We can take $a = (\ln \pi)/2$ to satisfy (\ref{EqXFPPOneHolds}) for all $x > a$. Observe that $\lim_{x \longrightarrow \infty} a_F(x) = (\ln \pi)/2$ by (\ref{EqLimitAF}).

$F(x) = e^\frac{1}{(12x+1)x}$. Here $A = C = 1$ and $B = 12$. We can take $a = 1/\sqrt{6}$. $\lim_{x \longrightarrow \infty} a_F(x) = 0$.

$F(x) = e^\frac{1}{12x^2}$. Here $A = 1$, $B = 12$, and $C = 0$. We can take $a = 1/\sqrt{6}$. $\lim_{x \longrightarrow \infty} a_F(x) = 0$.

Note that $L(x)$ is the product of the functions $F(x)$ in the first and second cases and $R(x)$ is the product of them in the first and third. Suppose $a(x), b(x)$ are positive valued strictly decreasing functions on $(a, \infty)$. Since sums and products of positive valued strictly decreasing functions are the same, the function $c(x) = a(x) + b(x) + a(x)b(x)/x$ is positive valued and strictly decreasing. Observe that $(1 + a(x)/x)(1+ b(x)/x) = 1 + c(x)/x$ for all $x > a$. Since $1/\sqrt{6} < (\ln \pi)/2$:
\begin{Prop}\label{PropLRAOverX}
Let $F(x) = L(x)$ or $F(x) = R(x)$. Then:
\begin{enumerate}
\item[{\rm (a)}] $a_F(x)$ is a strictly decreasing function on $((\ln \pi)/2, \infty)$.
\item[{\rm (b)}] $\lim_{x \longrightarrow \infty}a_F(x) = (\ln \pi)/2$.
\end{enumerate} \qed
\end{Prop}
\medskip

Observe that $0.5723 < (\ln \pi)/2 < 0.5724$. Estimates of $R(n)$ for $n = 1$, $10$, and $100$ are: $1.92648 < R(1) < 1.92649$, $1.05978 < R(10) < 1.05979$, and finally $1.005748 < R(100) < 1.005749$.  By Proposition \ref{PropLRAOverX} we have:
\begin{Cor}\label{LemmaRnOverUnderEstimate}
$1 + (\ln \pi)/(2x) < R(x)$ for $x \geq (\ln \pi)/2$ and $R(x) < 1 + a/x$ for all $x \geq b$, where $a = .9265$ and $b = 1$, $a = .5979$ and $b = 10$, or $a = .5749$ and $b = 100$.
\qed
\end{Cor}
\medskip

We leave the analysis of $P(x)$ to the reader. The function $a_P(x)$ is strictly increasing on $[1, \infty)$ and $\lim_{x \longrightarrow \infty}a_P(x) = \infty$. Strictly increasing comes down to $P(x) > \DS{\frac{x}{x + 1 - \ln x}}$ for $x > 1$ which can be shown in various ways. Estimates for $a_P(x)$ can be obtained using Lemma \ref{LemmaRNOverLN} or Proposition \ref{PropXToOneOverX}.

Finally, let $F(x) = T'(x) = \DS{P(x)\left(\frac{x+1 - \ln x}{x}\right)}$, where $a = e$. The properties of $\textsl{a}_s(x) = a_{T'}(x)$ we need in Section \ref{SecTDifffences} are easily developed there. Here we outline a deeper analysis of $\textsl{a}_s(x)$ analysis for the reader who might be interested. First of all
\begin{equation}\label{EqATPrime}
\textsl{a}_s(x) = P(x)(x + 1 - \ln x) - x
\end{equation}
and, since $P'(x) = p(x)P(x)$, where $p(x) = (1 - \ln x)/x^2$, it follows that
\begin{equation}\label{EqATPrimeTPrime}
\textsl{a}_s'(x) = \DS{\left(\frac{x^2 + (\ln x - 1)^2 - x\ln x}{x^2}\right)}P(x) - 1.
\end{equation}
As $x > \ln x$ for all $x > 0$ in fact, it follows that the coefficient of $P(x)$ in the preceding equation is positive.

We analyze the derivative $\textsl{a}_s'(x)$ in order to understand $\textsl{a}_s(x)$. Let
\begin{eqnarray*}
G(x) & = & \ln P(x) - \DS{\ln \left(\frac{x^2}{x^2 + (\ln x - 1)^2 - x\ln x}\right)} \\
& = & \DS{\frac{\ln x}{x} - 2\ln x + \ln (x^2 + (\ln x - 1)^2 - x\ln x)}.
\end{eqnarray*}
Then $\textsl{a}_s'(x) = 0$ if and only if $G(x) = 0$. The sign of $\textsl{a}_s'(x)$ is that of $G(x)$.

Let $y = \ln x$. Then $x = e^y$ and $y > 1$. Let $K(y) = \DS{e^y\left(\frac{4 - y}{(y-1)^2}\right)}$. Then $K(y)$ is a strictly decreasing function on its domain $(1, \infty)$ since its derivative is always negative and its range it the set of real numbers. Now
$$
G'(x) = \DS{\left(\frac{\ln x - 1}{x^2}\right)\left(\frac{(K(y)-1)(y-1)^2}{x^2 + (\ln x - 1)^2 - x\ln x}\right)}.
$$
Thus $G'(x)$ has a unique zero $x_0$, given by $K(\ln x_0) = 1$. Further $G'(x) > 0$ for $e < x < x_0$ and $G'(x) < 0$ for $x > x_0$. As $\lim_{x \longrightarrow \infty} G(x) = 0$ it follows that $G(x) > 0$ for $x \geq x_0$. Now $G(25.8679) < 0 < G(25.8680)$ shows that $G(x)$ has a unique zero $z$ and $25.8679 < z < 25.8680$. Furthermore $G(x) < 0$ for $e < x < z$ and $G(x) > 0$ for $x > z$. Thus $\textsl{a}_s(x)$ is strictly decreasing on $(e, z]$ and is strictly increasing on $[z, \infty)$. Hence $\textsl{a}_s(z)$ is an absolute minimum for $\textsl{a}_s(x)$.

Observe that $0.9114 < \textsl{a}_s(25.8679), \textsl{a}_s(25.8680), \textsl{a}_s(26) < 0.9115$. Thus $\textsl{a}_s(n) > 0.9114$ for all $n > e$, indeed for all $n \geq 1$.

We note $45.8750 < x_0 < 45.8751$ as $K(\ln 45.8750) > 1 > K(\ln 45.8751)$.
%
%SECTION
%
\section{The differences $s_n = (n+1)\sqrt[n+1]{n+1} - n\sqrt[n]{n}$}\label{SecTDifffences}
Let $n \geq 1$. We shall see that sequence of differences $S_n = e\sqrt[n+1]{(n+1)!} -  e\sqrt[n]{n!}$ is closely related to the relationship between powers and factorials in Section \ref{SecFactorialsandPowers}. We examine this sequence in Section \ref{SecTauDifffences}. In this section we consider the related sequence of differences $s_n = T(n+1) - T(n) = (n+1)\sqrt[n+1]{n+1} - n\sqrt[n]{n}$. See (\ref{EqCapSN}).

Let $x > 0$. Recall $\textsl{a}_s(x) = (T'(x) - 1)x$ and that $\DS{T'(x) = 1 + \frac{\textsl{a}_s(x)}{x}}$. Since $\textsl{a}_s'(x) = T''(x)x + T'(x) - 1$ it follows by (\ref{EqTPrime}) and (\ref{EqTPrimePrime}) that
$$
\textsl{a}_s'(x) = P(x)\left(\left(\frac{\ln x - 1}{x}\right)^2 - \frac{\ln x}{x} + 1\right) - 1
$$
and therefore
\begin{equation}\label{EqAsPrimeLimit}
\lim_{x \rTo \infty}\textsl{a}_s'(x) = 0
\end{equation}
as $\lim_{x \rTo \infty} P(x) = 1$ and $\DS{\lim_{x \rTo \infty}\frac{\ln x}{x} = 0}$. We have shown in Section \ref{SecERootNFactorial} that $T'(x)$ is strictly decreasing and $T'(x) > 1$ on $[1, \infty)$ and $\lim_{x \longrightarrow \infty}T'(x) = 1$. Writing $\DS{\textsl{a}_s(x) = \frac{T'(x) - 1}{(1/x)}}$, it follows by L'H\^{o}pital's rule and (\ref{EqTPrimePrime}) that
\begin{eqnarray*}
\lim_{x \rTo \infty}\textsl{a}_s(x) & = & \lim_{x \rTo \infty}-x^2T''(x)\\
& = & \lim_{x \rTo \infty} -P(x)\left(\frac{(\ln x - 1)^2}{x} - 1\right) = (-1)1(-1) = 1.
\end{eqnarray*}

Now $s_n = T'(x)$ for some $n < x < n+1$ by the Mean Value Theorem. In particular we have $1 < T'(n+1) < s_n < T'(n)$, that is
\begin{equation}\label{EqnSNEstimates}
1 < \sqrt[n+1]{n+1}\left(\frac{n+2-\ln (n+1)}{n+1}\right) < s_n < \sqrt[n]{n}\left(\frac{n+1-\ln n}{n}\right).
\end{equation}
\begin{Theorem}\label{TheoremNRootN}
Let $s_n = (n+1)\sqrt[n+1]{n+1} - n\sqrt[n]{n}$ for $n \geq 1$. Then:
\begin{enumerate}
\item[{\rm (a)}] $\lim_{n \longrightarrow \infty} s_n = 1$.
\item[{\rm (b)}] $2 > s_1 > s_2 > s_3 > \cdots > 1$.
\item[{\rm (c)}] $\DS{1 + \frac{\textsl{a}_s(n+1)}{n+1} < s_n} < 1 + \frac{\textsl{a}_s(n)}{n}\;\;$ for all $n \geq 1$.
\item[{\rm (d)}] $\lim_{n \longrightarrow \infty} \textsl{a}_s(n) = 1$.
\end{enumerate}
\end{Theorem}

\pf
Parts (a) and (b) follow by (\ref{EqnSNEstimates}), its preceding statement, and $\lim_{x \rTo \infty}T'(x) = 1$. Part (c) a restatement of of the last two inequalities of (\ref{EqnSNEstimates}). Part (d) follows from the limit calculation preceding (\ref{EqnSNEstimates}).
\qed
\medskip

Note that $\textsl{a}_s(x) = P(x)(x+1 - \ln x)$ by (\ref{EqTPrime}). Therefore
\begin{equation}\label{EqAXSubSLessOne}
\textsl{a}_s(x) < 1
\end{equation}
for all $x \geq 6.7537$ by Proposition \ref{PropXToOneOverX}. We will need
\begin{equation}\label{EqANSubSLessOne}
\textsl{a}_s(n) < 1
\end{equation}
for all $n \geq 1$. This inequality holds for $n \geq 7$ by (\ref{EqAXSubSLessOne}).  One can can show that $\textsl{a}_s(1), \ldots, \textsl{a}_s(6) < 1$ by direct computation.
%
%SECTION
%
\section{The differences $S_n = e\sqrt[n+1\;\;]{(n+1)!} -  e\sqrt[n]{n!}$}\label{SecTauDifffences}
Here we list some basic properties of the sequence $S_1, S_2, S_3, \ldots \;$ introduced at the beginning of Section \ref{SecTDifffences}. We use results from \cite{Radford}. Note the parallels between this section and the preceding one.

Let $x> 0$ and let $\Gamma(x)$ denote the gamma function. We refer the reader to \cite{Artin} as a basic reference for the gamma function. Set $G(x) = \Gamma(x+1)^{1/x}$, $A(x) = eG'(x)$, and $\textsl{a}_S(x) = (A(x) - 1)x$. Again, $\DS{A(x) = 1 + \frac{\textsl{a}_S(x)}{x}}$. Let $n \geq 1$. Since $G(n) = \sqrt[n]{n!} \;$, it follows that  $S_n = eG(n+1) - eG(n) = A(x)$ for some $n < x < n+1$ by the Mean Value Theorem. Now by part (a) of \cite[Proposition 1]{Radford} it follows that $A(x)$ is strictly decreasing on $[18, \infty)$. Thus $A(n) > S_n > A(n+1)$ for $n \geq 18$.
\begin{Theorem}[\cite{Radford}, Theorem 2]\label{ThmSNMain}
Let $S_n = e\sqrt[n+1]{(n+1)!} - e\sqrt[n]{n!} \;$ for $n \geq 1$. Then:
\begin{enumerate}
\item[{\rm (a)}] $\lim_{n \rTo \infty}S_n = 1$.
\item[{\rm (b)}] $1.15 > S_1 > S_2 > S_3 > \cdots > 1$.
\item[{\rm (c)}] $1 + \DS{\frac{\textsl{a}_S(n+1)}{n+1}} < S_n < 1 + \DS{\frac{\textsl{a}_S(n)}{n}}$ for all $n \geq 18$.
\item[{\rm (d)}] $\lim_{n \rTo \infty}\textsl{a}_S(n) = \DS{\frac{1}{2}}$.
\end{enumerate} \qed
\end{Theorem}
\medskip

We remark that $\lim_{x \rTo \infty}\textsl{a}_S'(x) = 0$ which is \cite[(24)]{Radford}. We need a fixed numerical value for an over estimate of $S_n$ in place of $\textsl{a}_S(n)$ in part (c) of Theorem \ref{ThmSNMain}. Our desired over estimate is derived from
\begin{equation}\label{EqSNDifferenceOver}
S_n < (1/2)(R(n+1)T(2(n+1)) - L(n)T(2n))
\end{equation}
which follows by (\ref{EqCapSN}). Note that $2 < s_{2n+1} + s_{2n} < 2(1 + 1/(2n))$ by part (c) of Theorem \ref{TheoremNRootN} from which we calculate $2 < T(2(n+1))) - T(2n) < 2(1 + 1/(2n))$. Recall that $R(n+1) < L(n)$ by (\ref{EqLnRn}). In light of (\ref{EqSNDifferenceOver}) we have:
\begin{equation}\label{EqSnLnEstimate}
S_n < \frac{1}{2}R(n+1)(T(2(n+1)) - T(2n)) < R(n+1)\left(1 + \frac{1}{2n}\right).
\end{equation}
In Section \ref{SecsigmaSequence} we use
\begin{equation}\label{EqSnLessThan2}
\DS{S_n < 1 + \frac{1.1}{n}}
\end{equation}
To establish this, we first observe that $S_n < 1 + \DS{\frac{0.5 + a_R(n+1)}{n}}$ which holds by (\ref{EqSnLnEstimate}), the definition of $a_R(x)$, and the fact that $a_R(n) > 0$. As $0.5 + a_R(n+1)$ is a strictly decreasing function on $[1, \infty)$ by part (a) of Proposition \ref{PropLRAOverX}, the calculations preceding Corollary \ref{LemmaRnOverUnderEstimate} establish the over estimate of (\ref{EqSnLessThan2}) with $1.1$ holds for $n \geq 10$. One can compute $\textsl{a}_S(1), \ldots, \textsl{a}_S(9)$ to complete the proof of (\ref{EqSnLessThan2}).
\section{A sequence $\sigma_1 \leq \sigma_2 \leq \sigma_3 \leq \cdots \;$ of positive integers}\label{SecsigmaSequence}
Suppose that $T_1, T_2, T_3, \ldots$ is any sequence of real numbers such that $T_1 \geq 1$ and the sequence $s_1, s_2, s_3, \ldots \;$, where $s_n = T_{n+1} - T_n$, satisfies
{\flushleft (s.1) $1 \leq s_n < 2$}
\medskip

\noindent
for all $n \geq 1$. Then
\begin{equation}\label{EqLowerTnUpper}
n \leq T_n < T_1+ 2(n-1)
\end{equation}
for all $n > 1$ by (s.1); the first inequality holds for all $n \geq 1$. Observe that the sequence $T_1', T_2', T_3', \ldots$, whose terms are defined by $T_n' = T_n - d$ for all $n \geq 1$, where $d = T_1 - 1$, satisfies $T_1' = 1$ and $s_n = T_{n+1}' - T_n'$ for all $n \geq 1$.

The sequence $s_1, s_2, s_3, \ldots \;$ of Section \ref{SecTDifffences} satisfies (s.1), where $T_n = n\sqrt[n]{n}$, and $T_1 = 1$.  See part (b) of Theorem \ref{TheoremNRootN}. The sequence $S_1, S_2, S_3, \ldots \;$ of Section \ref{SecTauDifffences} satisfies (s.1), where $T_n = e\sqrt[n]{n!}$, and $T_1 = e$.  See part (b) of Theorem \ref{ThmSNMain}.

Let $n \geq 1$ and let $\sigma_n$ be the largest integer $\ell$ which satisfies $n + (\ell - 1) \leq T_n$. Then $\sigma_n \geq 1$, by $T_1 \geq 1$ and (\ref{EqLowerTnUpper}), and $n + \sigma_n - 1$ is the largest integer equal to or less than $T_n$. Since
$$
n + \sigma_n = n + (\sigma_n - 1) + 1  \leq T_n + 1 \leq T_{n+1}
$$
by (s.1) it follow that $n + \sigma_n \leq (n+1) + (\sigma_{n+1} - 1)$. Therefore
\begin{equation}\label{EqSigmaNNPlusOne}
\sigma_n \leq \sigma_{n+1}.
\end{equation}
Now $T_n < n + \sigma_n$ and $n + \sigma_{n+1} < T_{n+1}$; thus
$$
\sigma_{n+1} - \sigma_n = (n + \sigma_{n+1}) - (n + \sigma_n) < T_{n+1} - T_n < 2
$$
by (s.1) which implies
\begin{equation}\label{EqSigmaNNPlusOneAlt}
\sigma_{n+1} = \sigma_n \;\; \mbox{or} \;\; \sigma_{n+1} = \sigma_n +1
\end{equation}
by (\ref{EqSigmaNNPlusOne}). As a result, if $a$ is any real number then $\sigma_n < a$ implies $\sigma_{n+1} < a+1$.

Applications of interest to us involve the sequence $\sigma_{n-m}, \ldots, \sigma_n$ when $2 \leq m = \sigma_n < n$.
\begin{Lemma}\label{LemmaSigmaNLessN}
Suppose that $T_1, T_2, T_3, \ldots$ is a sequence of real numbers such that $T_1 = 1$ and (s.1) are satisfied. Then $\sigma_1 = 1$ and $\sigma_n < n$ for all $n > 1$.
\end{Lemma}

\pf
By definition $\sigma_1 = 1$ as $T_1 = 1$. Now $T_2 < 3$ by (\ref{EqLowerTnUpper}) and thus $\sigma_2 = 1$. By the remark following (\ref{EqSigmaNNPlusOneAlt}) we conclude that $\sigma_n < n$ for all $n > 1$.
\qed
\medskip

Suppose that $T_1, T_2, T_3, \ldots$ is a sequence of real numbers such that $T_1 \geq 1$ and (s.1) are satisfied. Let $n \geq 1$ and set $\nu_n = n + \sigma_n$. Then $\nu_n - 1 \leq T_n < \nu_\ell$,
$\nu_{n+1} = \nu_n + 1 + (\sigma_{n+1} - \sigma_n)$, and thus
\begin{equation}\label{EqNuSubEllPlusU}
\nu_{n+\ell} = \nu_n + \ell + (\sigma_{n + \ell} - \sigma_n)
\end{equation}
for all $\ell \geq 0$.
\begin{Lemma}\label{LemmaSigmasSame}
Let $T_1, T_2, T_3, \ldots$ be a sequence of real numbers such that $T_1 \geq 1$ and (s.1) are satisfied. Suppose that $2 \leq m = \sigma_n < n$ and $\sigma_{n-m} = \cdots = \sigma_n$. Then $n-1 \leq T_{n-m} < n \leq T_{n-m+1} < n+1 \leq T_{n-m+2} < n+2$.
\end{Lemma}

\pf
First of all $n-m+2 \leq n$ since $m \geq 2$. By assumption $n - m \geq 1$. Note that $\nu_{n-m+\ell} = \nu_{n-m} + \ell$ for all $0 \leq \ell \leq m$ by assumption and (\ref{EqNuSubEllPlusU}). In particular $\nu_n = \nu_{n-m} + m$ which means $\nu_{n-m} = n$. Therefore $\nu_{n-m+\ell} = n+\ell$ for all $0 \leq \ell \leq m$ and the lemma now follows. Or simply use the fact that $i + (\sigma_i - 1) \leq T_i < i + \sigma_i$ for all $i \geq 1$ to establish the lemma.
\qed
\medskip

The lemma is about the sequence $T_1, T_2, T_3, \ldots \;$ when $\sigma_{n-m}, \ldots, \sigma_n$, where $2 \leq m = \sigma_n < n$, has one value. In the two cases of interest to us, namely the sequences $s_1, s_2, s_3, \ldots \;$  and $S_1, S_2, S_3, \ldots \;$ of Sections \ref{SecTDifffences} and \ref{SecTauDifffences} respectively, the sequence $\sigma_{n-m}, \ldots, \sigma_n$, where $2 \leq m = \sigma_n < n$, has one or two values as we as show later. In our two cases $T_1, T_2, T_3, \ldots$ satisfies $T_1 \geq 1$, (s.1), and:
{\flushleft (s.2) There exists $a > 0$ such that $s_n < 1 + a/n$}
\medskip

\noindent
for all $n \geq 1$.  See (\ref{EqANSubSLessOne}) and (\ref{EqSnLessThan2}).

Let $T_1, T_2, T_3, \ldots$ satisfy $T_1 \geq 1$, (s.1), and (s.2). Then $\lim_{n \longrightarrow \infty} s_n = 1$ by (s.1) and (s.2).

Set $d_n = T_n - (n + \sigma_n - 1)$. Observe that $0 \leq d_n < 1$. Now
\begin{eqnarray*}
d_{n+1} - d_n & = & T_{n+1} - (n + \sigma_{n+1}) - (T_n - (n + \sigma_n - 1)) \\
  & = & T_{n+1} - T_n - 1 - (\sigma_{n+1} - \sigma_n) \\
  & < & a/n - (\sigma_{n+1} - \sigma_n) < a/n.
\end{eqnarray*}
The last two inequalities follow by (s.2) and (\ref{EqSigmaNNPlusOneAlt}) respectively.

Assume that $\sigma_{n+1} = \sigma_n + 1$. By the preceding calculation, $1 - a/n < d_n$, since $d_{n+1} \geq 0$, and $d_{n+1} < a/n$, since $d_n < 1$.

Let $r \geq 1$ and suppose, in addition, that $\sigma_{n + r + 1} = \sigma_{n+r} + 1$. Then $1 - a/(n+r) < d_{n+r}$ and $$d_{n+r} = \sum_{i = 1}^{r-1} (d_{n+i + 1} - d_{n+i}) + d_{n+1} < \sum_{i = 1}^{r-1} a/(n+i) + a/n < ra/n.$$ Therefore $1 - a/(n+r) < ar/n$, which can be rewritten $(1 + r/n)(n-ra) < a$ or $(1 + r/n)(n/a-r) < 1$. Consequently $n/a-r < 1$ which means $n/a < r+1$.
\begin{Prop}\label{PropVeryRare}
Let $T_1, T_2, T_3, \ldots$ be a sequence of real numbers such that $T_1 \geq 1$, (s.1), and (s.2) are satisfied. Let $n \geq 1$ and suppose $\sigma_{n+1} = \sigma_n + 1$. Then $\sigma_{n+1} = \sigma_{n+r}$ for all $r \geq 1$ or $\sigma_{n+1} = \cdots = \sigma_{n+r}$ and $\sigma_{n+r+1} = \sigma_{n+r} + 1$  for some $r \geq 1$. The inequality $r + 1 > n/a$ holds for any $r \geq 1$ which satisfies the last equation. \qed
\end{Prop}
\medskip

Let $T_1, T_2, T_3, \ldots$ be a sequence of real numbers such that $T_1 \geq 1$, (s.1), and (s.2) are satisfied. Let $n \geq 1$. Then $\sigma_n = \sigma_1 + \ell$ for some $\ell \geq 0$. Suppose $\ell \geq 1$. Then there is a sequence $1 \leq n_1 < n_2 < \cdots < n_\ell < n$, where $n_i$ is the largest integer $u$ which satisfies $\sigma_u = \sigma_1 + (i-1)$. For $1 \leq i < \ell$ let $r_i = n_{i+1} - n_i$. Then $\sigma_{n_i+1} = \cdots = \sigma_{n_i+r_i}$ and $r_i + 1 > n_i/a$ by Proposition \ref{PropVeryRare}. We have shown that $n_{i+1} > n_i(1 + 1/a) - 1$ for all $1 \leq i < \ell$. Therefore
$$
n_i \geq (1 + 1/a)^{i-1}(n_1 - a) + a
$$
for all $1 \leq i \leq \ell$, with strict inequality when $1 < i$. Since $n > n_\ell$ and $\sigma_n = \sigma_1 + \ell$:
\begin{Cor}\label{CorVeryRare}
Let $T_1, T_2, T_3, \ldots$ be any sequence of real numbers such that $T_1 \geq 1$, (s.1), and (s.2) are satisfied. Suppose $n_1 \geq 1$ satisfies $\sigma_{n_1} = \sigma_1$ and $\sigma_{n_1+1} = \sigma_1 + 1$. Then $n > (1 + 1/a)^{ \sigma_n - \sigma_1-1}(n_1 - a) + a$ for all $n > n_1$. \qed
\end{Cor}
\medskip

Suppose the sequence $\sigma_{n-m}, \ldots, \sigma_n$, where $2 \leq m = \sigma_n < n$, has more than two values. We appeal to Corollary \ref{CorVeryRare}. There are integers $m_1$ and $m_2$ such that $n-m \leq m_1 < m_2 < n$, $\sigma_{n-m} =  \cdots = \sigma_{m_1}$, $\sigma_{m_1+1} = \sigma_{m_1} + 1 = \cdots = \sigma_{m_2}$, and $\sigma_{m_2+1} = \sigma_{m_2} + 1 \leq \sigma_n$. In particular there is an $n_1 \geq 1$ such that $\sigma_{n_1} = \sigma_1$ and $\sigma_{n_1+1} = \sigma_1 + 1$. Observe that $n_1 \leq m_1 < n$.

Let $r = m_2 - m_1$. Then $r \geq 1$ and $\sigma_{m_1 + r + 1} = \sigma_{m_1+r}+1$. Therefore $m_2 - m_1 + 1 > m_1/a$ by Proposition \ref{PropVeryRare}. We have shown that
$$
n \geq m_2 + 1 > m_1(1+1/a) \geq (n-m)(1+1/a)
$$
and therefore
\begin{equation}\label{EqVeryRare}
m(1+a) > n > (1+1/a)^{m-\sigma_1-1}(n_1-a)+a
\end{equation}
by Corollary \ref{CorVeryRare}.

We specialize to the case $T_n = n\sqrt[n]{n}$ for $n \geq 1$ to which the preceding corollary applies. Suppose $\sigma_{n-m}, \ldots, \sigma_n$, where $2 \leq m = \sigma_n < n$, has more than two values. Then (\ref{EqVeryRare}) holds which we use to arrive at a contradiction.  We can take $a = 1$ for (s.2) by (\ref{EqANSubSLessOne}). Observe $\sigma_1 = \sigma_2 = 1$, $\sigma_3 = \sigma_4 = \sigma_5 = 2$, $\sigma_6 = \cdots = \sigma_{15} = 3$, and $\sigma _{16} = 4$. In particular $n_1 = 2$ and $2m > n > 2^{m-2}+1$. As a result $2m \geq 2^{m-2} + 3$ which implies $2 \leq m \leq 4$. This means $2 < n \leq 7$. Since $\sigma_7 = 3$, we have shown $m = 2$ or $m = 3$ and consequently $n \leq 5$. In this case the values of the $\sigma_{n-m}, \ldots, \sigma_n$ are among $1$ and $2$, a contradiction. Thus:
\begin{Cor}\label{CorOneOrTwoValuesA}
Let $T_1, T_2, T_3, \ldots$ be the sequence whose terms are defined by $T_n = n\sqrt[n]{n}$ for $n \geq 1$ for all $n \geq 1$. If $n \geq 3$ then $2 \leq m = \sigma_n < n$ and the sequence $\sigma_{n-m}, \ldots, \sigma_n$ has one or two values. \qed
\end{Cor}
\medskip

We continue with the sequence $T_1, T_2, T_3, \ldots$ of the preceding corollary. Observe that the increasing sequence $\sigma_1, \sigma_2, \sigma_3, \ldots$ is unbounded. Suppose not. Then there is an $M > 0$ such that $T_n < n + M$ for all $n \geq 1$ and thus $\ln n < a_P(n) < M$ all $n > 1$, contradiction. See the paragraph which precedes Proposition \ref{PropXToOneOverX}. We have shown that $\sigma_1, \sigma_2, \sigma_3, \ldots$ is unbounded.

Let $n_0 = 0$ and for $i \geq 1$ let $n_i$ be the largest integer $u$ which satisfies $\sigma_u = \sigma_1 + (i - 1)$. Then $n_0, n_1, n_2, n_3, \ldots$ is a strictly increasing sequence which is determined by $\sigma_{n_i + 1} = \cdots = \sigma_{n_{i+1}}$ for all $i \geq 0$ and $\sigma_{n_i + 1} = \sigma_{n_i}+1$ for all $i \geq 1$. Then $n_{i+1} - n_i + 1 > n_i/1$, or $n_{i+1} \geq 2n_i$, for all $i \geq 1$ by Proposition \ref{PropVeryRare}. As $n_1 = 2$ it follows that $n_i \geq 2^i$ for all $i \geq 1$. It is interesting to note that $n_1 = 2, n_2 = 5$, and $n_3 = 15$.
\begin{Prop}\label{PropSigmasTwoValues}
Let $T_1, T_2, T_3, \ldots$ be a sequence of real numbers such that $T_1 \geq 1$, (s.1), and (s.2) are satisfied. Let $n \geq 1$ and suppose that the sequence $\sigma_{n-m}, \ldots, \sigma_n$ has two values, where $2 \leq m = \sigma_n < n$. Then:
\begin{enumerate}
\item[{\rm (a)}] There exists an $\ell$ which satisfies $n - m \leq \ell < n$, $\sigma_{n-m} = \cdots = \sigma_\ell = m-1$, and $\sigma_\ell + 1 = \sigma_{\ell + 1} = \cdots = \sigma_n = m$.
\item[{\rm (b)}] $n-2 \leq T_{n-m} < n-1$.
\item[{\rm (c)}] If $\ell = n - m$ then $n \leq T_{n-m+1} < n+1 \leq T_{n-m+2} < n+2$
\item[{\rm (d)}] If $\ell = n - m + 1$ then $n-1 \leq T_{n-m+1} < n$ and $n+1 \leq T_{n-m+2} < n+2$.
\item[{\rm (e)}] If $\ell = n - m + 2$ then $n-1 \leq T_{n-m+1} < n \leq T_{n-m+2} < n+1$ and $n+2 \leq T_{n-m+3} < n+3$.
\item[{\rm (f)}] $n-1 \leq T_{n-m+1} < n \leq T_{n-m+2} < n+1 \leq T_{n-m+3} < n+2$ when $\ell \geq n - m + 3$.
\end{enumerate}
\end{Prop}

\pf
Part (a) follows by (\ref{EqSigmaNNPlusOne}) and (\ref{EqSigmaNNPlusOneAlt}). Note that $n-m+2 \leq n$ since $m \geq 2$. Using the fact that $i + (\sigma_i - 1) \leq T_i < i + \sigma_i$ for all $i \geq 1$, parts (b)-(f) are easily established.
\qed
\medskip

Let $n \geq 1$. We specialize to the case $T_n = e\sqrt[n]{n!}$ to which Corollary \ref{CorVeryRare} applies. Recall that $R(n) < 1 + 1/n$ by Corollary \ref{LemmaRnOverUnderEstimate}. For ease of calculation, we may take $a = 3/2$ for (s.2) by (\ref{EqSnLessThan2}). Using (\ref{EqCapSN}) one can easily show that $\sigma_1 = \sigma_2 = \sigma_3 = 2$, $\sigma_4 = \cdots = \sigma_{54} = 3$, and $\sigma_{55} = 4$. In particular $2 \leq m = \sigma_n< n$ for all $n \geq 3$. Note that $\sigma_1 = 2$ and $n_1 = 3$.

Suppose $2 \leq m = \sigma_n < n$ and $\sigma_{n-m}, \ldots, \sigma_n$ has more than two values. Then $(5/2)m > n > (5/3)^{m - 3}(3/2)$ by (\ref{EqVeryRare}) and therefore $m > (5/3)^{m-4}$. This means $m \leq 8$ and therefore $n \leq 20$. If $n \leq 20$, the values of the sequence $\sigma_{n-m}, \ldots, \sigma_n$ are $2$ or $3$ as $\sigma_1 = 2$ and $\sigma_{20} = 3$, a contradiction. Thus in any case the sequence $\sigma_{n-m}, \ldots, \sigma_n$, where $2 \leq m = \sigma_n < n$, has one or two values.
\begin{Cor}\label{CorOneOrTwoValuesB}
Let $T_1, T_2, T_3, \ldots$ be the sequence whose terms are defined by $T_n = e\sqrt[n]{n!}$ for all $n \geq 1$. If $n \geq 3$ then $m = \sigma_n < n$ and the sequence $\sigma_{n-m}, \ldots, \sigma_n$ has one or two values.
 \qed
\end{Cor}
\medskip

We continue with the sequence $T_1, T_2, T_3, \ldots$ of the preceding corollary. Observe that the increasing sequence $\sigma_1, \sigma_2, \sigma_3, \ldots$ is unbounded. Suppose this is not the case. Then there is an $M > 0$ such that $T_n < n + M$ for all $n \geq 1$. Since $L(x) > 1$ for all $x > 0$, we can use (\ref{EqCapSN}) and the discussion preceding Proposition \ref{PropXToOneOverX} to conclude $\ln 2n < a_P(2n) < 2M$  for all $n > 1$. The latter is not the case for large $n$. Therefore $\sigma_1, \sigma_2, \sigma_3, \ldots$ is unbounded.

Let $n_1, n_2, n_3, \ldots$ be the strictly increasing sequence derived from the sequence $\sigma_1, \sigma_2, \sigma_3, \ldots$ as was its analog in the discussion after Corollary \ref{CorOneOrTwoValuesA}. Then $n_{i+1} - n_i + 1 > n_i/(3/2)$, or $n_{i+1} > (5/3)n_i - 1$, for all $i \geq 1$ by Proposition \ref{PropVeryRare}. From $3n_{i+1} > 5n_i - 3$ we conclude that $3n_{i+1} \geq 5n_i - 2$ and thus $n_{i+1} \geq (5/3)n_i - 2/3$ for all $i \geq 1$. By induction $n_i \geq 2(5/3)^{i-1} + 1$ for all $i \geq 1$ follows. It is interesting to note that $n_1 = 3$ and $n_2 = 54$.

\section{Factorials and powers}\label{SecFactorialsandPowers}
Let $a > 1$ and $n$ be a positive integer. We examine the relationship between $n!$ and $a^n$. First of all observe that $n < a$ implies $n! < a^n$, or equivalently $n! \geq a^n$ implies $n \geq a$.

Fix $a$. There exists an $n$ such that $n! \geq a^n$; consider the power series expansion of $e^a$ for example. Let $n_a$ be the smallest such positive integer. Then $n_a \geq 2$ since $a > 1$. Note that $n_a$ is an integer $n$ which satisfies $n \geq 2$, $a^n \leq n!$, and $(n-1)! < a^{n - 1}$, or equivalently
\begin{equation}\label{EqNSubA}
n \geq 2 \;\;\mbox{and} \;\; \sqrt[n-1]{(n-1)!} < a \leq \sqrt[n]{n!}.
\end{equation}
At this point we consider the sequence $1 = t_1, t_2, t_3, \ldots \;$, where $t_n = \sqrt[n]{n!}\,$. By part (b) of Theorem \ref{ThmSNMain}, the sequence is strictly increasing and unbounded. There are simpler arguments for these assertions. Unbounded follows by (\ref{EqNSubA}) with $n = n_a$.  The calculation
$$
(t_{n+1}/t_n)^{n+1} = (n+1)/\sqrt[n]{n!} \geq (n+1)/\sqrt[n]{n^n} = 1 + 1/n > 1
$$
shows that the sequence is strictly increasing. Now (\ref{EqNSubA}) is equivalent to $n \geq 2$ and $a \in (t_{n-1}, t_n]$. Therefore:
\begin{Lemma}\label{LemmaNsubA}
Let $a > 1$. Then (\ref{EqNSubA}) has a unique solution $n$ which is $n_a$.
\qed
\end{Lemma}
\medskip

Observe that $n! < a^n$ when $n < n_a$, $a^n \leq n!$ when $n = n_a$, and $a^n < n!$ when $n_a < n$. Let $n \geq 2$ and set $a = \sqrt[n]{n!}$\,. Then $a > 1$ and $n = n_a$ by our last assertion. This means the $n_a$'s constitute the positive integers greater than one. Observe that $a < b$ implies $n_a \leq n_b$. As a result $\lim_{a \longrightarrow \infty} n_a = \infty$.

The values $t_n$ are closely related to the classical approximations of $e$. We recall that $(1 + 1/n)^n < e < (1 + 1/n)^{n+1}$ for all $n \geq 1$. The terms $(1 + 1/n)^n$ form a strictly increasing sequence which converges to $e$ and the terms $(1 + 1/n)^{n+1}$ form a strictly decreasing sequence which converges to $e$. For $n \geq 1$ set $\epsilon_n = (1 + 1/n)^n/e$ and $\delta_n = \epsilon_1 \cdots \epsilon_n$. Then $0 < \epsilon_n, \delta_n < 1$, and
\begin{equation}\label{EqDeltaN}
\DS{\delta_n = \frac{(n+1)^n}{n!e^n}}
\end{equation}
follows by induction on $n$. Therefore
\begin{equation}\label{EqNRootNFact}
\DS{\frac{n+1}{e} = \sqrt[n]{\delta_nn!} < \sqrt[n]{n!} = \frac{n+1}{e\sqrt[n]{\delta_n}},}
\end{equation}
in particular $n+1 < e\sqrt[n]{n!}$, for all $n \geq 1$.
 The formula
\begin{equation}\label{EqRootNFact2}
\DS{\lim_{n \longrightarrow \infty} \frac{n}{\sqrt[n]{n!}} = e}
\end{equation}
follows from (\ref{EqCapSN}) and $\lim_{x \longrightarrow \infty} L(x) =  \lim_{x \longrightarrow \infty} P(x) = \lim_{x \longrightarrow \infty} R(x) = 1$. From (\ref{EqDeltaN}) and (\ref{EqRootNFact2}) we deduce
\begin{equation}\label{EqLimRootDeltaN}
\lim_{n \longrightarrow \infty} \sqrt[n]{\delta_n} = 1.
\end{equation}
In light of (\ref{EqNRootNFact}) and (\ref{EqLimRootDeltaN}) note the asymptotic relationship $\sqrt[n]{n!} \sim n/e$. This makes a connection between the set of half open intervals $(\sqrt[n]{n!}, \sqrt[n+1]{(n+1)!}]$ implicit in (\ref{EqNSubA}) and the set of half open intervals $(n/e, (n+1)/e]$ arising in Theorem \ref{TheoremMain} below.

Note by part (a) of Theorem \ref{ThmSNMain} the sequence of lengths of the first set of intervals converges to $1/e$ which is the length of each of the intervals in the second set. By part (b) of the same, the length of each interval in the first set exceeds $1/e$.

Using Lemma \ref{LemmaNsubA} and (\ref{EqRootNFact2}) one can show that

\begin{equation}\label{EqNaOverA}
\lim_{a \longrightarrow \infty}\frac{n_a}{a} = e.
\end{equation}
The main result of this paper is:
\begin{Theorem}\label{TheoremMain}
Suppose $n \geq 3$ and $n/e < a \leq (n+1)/e$. Let $m$ be the integer determined by $n + m - 1 \leq e\sqrt[n]{n!} < m + n$. Then:
\begin{enumerate}
\item[{\rm (a)}] $2 \leq m = \sigma_n < n$ and the segment $\sigma_{n-m}, \ldots, \sigma_n$ has one or two values.
\end{enumerate}
When the segment has two values let $n-m \leq \ell < n$ be the solution to $\sigma_{\ell + 1} = \sigma_\ell +1$.
\begin{enumerate}
\item[{\rm (b)}] If the segment has one value or has two values and $\ell = n-m$ then $n_a = n-m+1$ or $n_a = n-m+2$.
\item[{\rm (c)}] If the segment has two values and $\ell = n-m+1$ then $n_a = n-m+2$.
\item[{\rm (d)}] If the segment has two values and $\ell \geq n-m+2$ then $n_a = n-m+2$ or $n_a = n-m+3$.
\end{enumerate}
\end{Theorem}

\pf
First note that $a > 1$. Let $T_n = e\sqrt[n]{n!} = et_n$ for all $n \geq 1$. Then $T_1, T_2, T_3, \ldots$ is a strictly increasing sequence since $t_1, t_2, t_3, \ldots$ is. We have noted that $T_1, T_2, T_3, \cdots \;$ satisfies $T_1 \geq 1$, (s.1), and (s.2). Let $n \geq 3$. We established in the discussion following Proposition \ref{PropSigmasTwoValues} that $2 \leq m =\sigma_n < n$. By Corollary \ref{CorOneOrTwoValuesB} the sequence $\sigma_{n-m}, \ldots, \sigma_n$ has one or two values.

Suppose $\sigma_{n-m}, \ldots, \sigma_n$ has one value. Then $T_{n-m} < n$ and $n+1 \leq T_{n-m+2}$ by Lemma \ref{LemmaSigmasSame}. Therefore $\sqrt[n-m]{(n-m)!} < a \leq \sqrt[n-m+2]{(n-m+2)!}$ which means $n_a = n-m+1$ or $n-m+2$ by (\ref{EqNSubA}), a result we use implicitly in the remainder of the proof.

Suppose $\sigma_{n-m}, \ldots, \sigma_n$ has more than one value. Then it has two values and Proposition \ref{PropSigmasTwoValues} applies. Let $\ell$ be as in part (a) this proposition.

First of all, suppose that $\ell = n-m$. Then by parts (b) and (c) of Proposition \ref{PropSigmasTwoValues} we have $T_{n-m} < n-1$ and $n+1 \leq T_{n-m+2}$  respectively. As a result $\sqrt[n-m]{(n-m)!} < a \leq \sqrt[n-m+2]{(n-m+2)!}$ and therefore $n_a = n-m+1$ or $n-m+2$ again.

Suppose $\ell = n-m+1$. Then $T_{n-m+1} < n$ and $n+1 \leq T_{n-m+2}$  by part (d) of the same. Therefore $\sqrt[n-m+1]{(n-m+1)!} < a \leq \sqrt[n-m+2]{(n-m+2)!}$ which means that $n_a = n-m+2$.

Suppose $\ell = n-m+2$. Then $T_{n-m+1} < n$ and $n+2 \leq T_{n-m+3}$ by part (e) of the same. Therefore $\sqrt[n-m+1]{(n-m+1)!} < a < \sqrt[n-m+3]{(n-m+3)!}$ which implies $n_a = n-m+2$ or $n-m+3$.

Suppose $\ell \geq n-m+3$. Then $T_{n-m+1} < n$ and $n+1 \leq T_{n-m+3}$ by part (f) of the same. Therefore $\sqrt[n-m+1]{(n-m+1)!} < a \leq \sqrt[n-m+3]{(n-m+3)!}$ which implies $n_a = n-m+2$ or $n-m+3$ again.
\qed
\medskip

The situation not covered by the preceding theorem is $1 < a \leq 3/e$. Here $a^2 < 2$; therefore $n_a = 2$.

Apropos of the preceding theorem, an example where $n$ is very large. Let $n = 10^{12}$. Using (\ref{EqCapSN}) one can show that $n + 14 < e\sqrt[n]{n!}  < n + 15$, hence $m = \sigma_n = 15$, and
$$
e\sqrt[n-15]{(n-15)!} < n < e\sqrt[n-14]{(n-14)!} < n+1 < e\sqrt[n-13]{(n-13)!}.
$$
Therefore $n_a = n-14$ or $n-13$, where $n/e < a \leq (n+1)/e$. In either case $n_a$ differs from $n$ by no more than $0.0000000014\%$. Perhaps the reader will find the following table intriguing.

\begin{center}
\begin{tabular}{|l|l|l|l|l|l|l|l|r|r|l|l|l|} \hline
$n$ & 1 & 2 & 3 & 4 & 5 & 6 & 7 & 8 & 9 & 10& 11 & 12 \\ \hline
$\sigma_{10^n}$  & 3 & 4 & 5 & 6 & 7 & 8 & 9 & 11  & 12 & 13 & 14 & 15 \\  \hline
\end{tabular}
\end{center}
\medskip

We close with comments on the relationship between half open intervals of the type $(t_n, t_{n+1}]$, where $n \geq 1$, and those of the type $(n/e, (n+1)/e]$, where $n \geq 3$. We have noted that the length of an interval of the first type exceeds the length of any interval of the second type. Thus an interval of the second type is contained in an interval of the first type or the union of two consecutive intervals of the first type. The proof of Theorem \ref{TheoremMain} makes this relationship explicit.


\begin{thebibliography}{99}
\bibitem{Artin}
Artin, Emil (translator Micheal Butler)
(2018).
The gamma function,
\emph{Dover Publications, Mineola, NY, USA}
%\textbf{43},
%pp. 225–-249.

\bibitem{MATLAB}
Mathworks, Inc. %Natick, MA,
(2020).
MATLAB 2020a
\emph{www.mathworks.com}
%\textbf{171},
%pp. 375–-417.

\bibitem{Radford}
Radford, David E.
(2021).
The gamma function and a certain sequence of differences,
\emph{preprint}
%\textbf{43},
%pp. 225–-249.

\bibitem{Robbins}
Robbins, Herbert
(1955).
A remark on Sterling's formula,
\emph{The American Mathematical Monthly}
\textbf{62},
pp. 26--29.

\bibitem{Wells}
Wells, Christian Nathan
(2017).
Private communication.
\end{thebibliography}
\end{document}